\documentclass[11pt,french]{amsart}
\usepackage[francais]{babel}
\usepackage[utf8]{inputenc}
\usepackage{latexsym}
\usepackage[all]{xy}
\usepackage{amsmath}
\usepackage{amssymb}
\usepackage{amsfonts}

\newcommand {\OO} {\mathcal{O}}

\newcommand {\SH} {\mathcal{SH}}
\newcommand {\D} {\mathcal{D}}

\newcommand {\BU} {B\mathbb{U}}

\newcommand {\Sch} {\mathbf{Sch}}

\newcommand{\s}{\infty}

\newcommand {\Ql}{\mathbb{Q}_{\ell}}

\newtheorem{thm}{Th\'eor\`eme}[section]
\newtheorem{prop}[thm]{Proposition}

\newtheorem{df}[thm]{D\'efinition}

\newtheorem{conj}[thm]{Conjecture}

\begin{document}

\title{\textbf{Géométrie non-commutative, formule des traces et conducteur de Bloch}}
\bigskip
\bigskip


\author{Bertrand To\"en$^*$}
   \address{IMT UMR 5219, CNRS, Universit\'{e} Paul Sabatier, Toulouse - France}
   \email{bertrand.toen@math.univ-toulouse.fr}

\author{Gabriele Vezzosi}
   \address{Dipartimento di Matematica ed Informatica ``Ulisse Dini'', Universit\`a di Firenze, Firenze - Italia}
   \email{gabriele.vezzosi@unifi.it}
   
   \thanks{$*$ Partially supported by 
ANR-11-LABX-0040-CIMI within the program  ANR-11-IDEX-0002-02}

\bigskip 

\date{Septembre 2016}

\maketitle

\begin{abstract}
Ce texte est basé sur l'exposé du premier auteur au premier congrès de la SMF (Tours, 2016). 
On y présente la formule du conducteur de Bloch, qui est une formule conjecturale décrivant 
le changement de topologie dans une famille de variétés algébriques lorsque le paramètre 
se spécialise en une valeur critique. L'objectif de ce texte est de présenter une approche
générale à la résolution de cette conjecture basée sur des techniques de géométrie 
non-commutative et de géométrie dérivée. 
\end{abstract}

\tableofcontents

\section*{Introduction}

Ce texte est un survol qui présente un approche à la conjecture du conducteur de Bloch (voir
conjecture \ref{cb}) basée sur des méthodes et des idées de géométrie non-commutative et de géométrie dérivée. Cette 
conjecture, qui prédit la variation de la caractéristique d'Euler d'une famille
de variétés algébriques lorsque le paramètre tend vers une valeur critique, est de nature
géométrico-arithmétique, et les mathématiques de ce texte s'insèrent dans un
programme plus vaste qui consiste à utiliser les méthodes de la géométrie non-commutative
(voir notre $\S$ 2), mais aussi de la géométrie dérivée (voir par exemple \cite{ems}), 
pour approcher des questions classiques de géométrie arithmétique. Ce texte de survol, 
ainsi que les textes plus détaillés \cite{mf,trace}, en constituent le premier pas. \\

Dans \cite{bl}, Bloch introduit une formule appelée aujourd'hui \emph{formule 
du conducteur de Bloch}, et en donne une preuve en dimension relative 1 (pour une famille
de courbes algébriques). Le contexte général est celui d'une famille de variétés
algébriques projectives $\{X_t\}_t$, 
qui varient algébriquement en fonction du paramètre $t$, et que
l'on étudie localement autour d'une valeur critique $t=o$. On suppose que
les variétés $X_t$ sont toutes non-singulières. Lorsque $X_o$ est elle-même non-singulière
on sait que les cohomologies de $X_t$ et de $X_o$ sont isomorphes et leurs caractéristiques
d'Euler sont donc égales. Lorsque $X_o$ est éventuellement singulière cette caractéristique
d'Euler change et la formule du conducteur de Bloch est une formule qui exprime 
leur différence $\chi(X_o)-\chi(X_t)$ en termes géométrico-algébriques (voir
conjecture \ref{cb}). Il faut rajouter ici que l'on ne considère aucune restriction 
sur les corps de base des variétés algébriques, les variétés $X_t$ et $X_o$ peuvent par
exemple toutes deux êtres des variétés algébriques complexes, toutes deux définies sur
un corps $k$ de caractéristique positive, mais on peut aussi avoir des situations
où $X_t$ est de caractéristique nulle alors que $X_o$ est de caractéristique positive. 
La formule du conducteur de Bloch n'est ainsi pas seulement une formule d'origine géométrique
mais rend compte aussi de phénomènes arithmétiques liés à la notion de ramification 
sauvage. 

La formule du conducteur de Bloch est un théorème dans certains cas (rappelés
dans notre $\S$ 1), mais reste ouverte en général. 
Récemment, dans \cite{trace} nous avons proposé de rapprocher cette formule d'une formule
des traces dans le cadre de la géométrie non-commutative. Le but de ce texte de survol
est de présenter de manière plus abordable le contenu de \cite{trace}. L'idée générale
de cette approche est d'introduire un schéma non-commutatif associé à 
la famille $\{X_t\}_t$ dont la caractéristique d'Euler est exactement la différence
que l'on cherche à quantifier $\chi(X_o)-\chi(X_t)$. Une fois cet objectif atteint, on 
invoque une formule du type Gauss-Bonnet pour ce schéma non-commutatif afin 
de calculer sa caractéristique d'Euler en termes géométrico-algébrique. C'est cette idée 
générale que nous allons tenter de décrire dans ce texte. Nous commencerons 
par rappeler la conjecture du conducteur de Bloch, sa signification dans certains cas
particuliers et le cas connus. Dans la seconde section nous rappellerons une approche
à la notion de schéma non-commutatif basée sur les dg-catégories. Nous tentons d'y donner
les définitions nécessaires à la compréhension de la suite du texte mais évitons 
très largement les détails techniques d'algèbre homotopique et de théorie
des $\s$-catégories qui sont malheureusement indispensables. Dans le $\S 3$ nous
introduisons la cohomologie des schémas non-commutatifs. Une première partie concerne
la cohomologie de Hodge et est relativement classique et connue sous le nom d'homologie
de Hochschild. Dans une seconde partie nous présentons la cohomologie 
$\ell$-adique des schémas non-commutatifs, qui a été introduite récemment dans
\cite{mf} et dont la construction est plus sophistiquée. La section $\S$ 4 concerne la
formule des traces (de type Lefschetz) pour les schémas non-commutatifs, qui permet 
de calculer la trace d'un endomorphisme sur la cohomologie à l'aide
d'un nombre d'intersection qui se décrit comme une dimension d'homologie de Hochschild. 
Enfin, dans la dernière section nous présentons l'application de la formule des traces
à la conjecture du conducteur de Bloch. Nous expliquons en particulier
comment elle permet d'affirmer que la conjecture est vraie lorsque la monodromie
est unipotente. Nous donnons aussi quelques éléments pour convaincre que le cas
général doit pouvoir s'en déduire (mais les détails restent sous investigation actuellement). 
\\

Le contenu de ce texte est essentiellement celui de l'exposé du premier auteur
au \emph{Premier congrès de la Société Mathématiques de France}, qui s'est tenu
à Tours en Juin 2016. Nous remercions vivement le SMF pour l'organisation de cet événement.

\section{La formule du conducteur de Bloch}

La formule du conducteur de Bloch est une formule numérique qui exprime 
les changements topologiques intervenant dans une dégénérescence de variétés
algébriques. Elle a été introduite dans \cite{bl} pour une famille de courbes algébriques, 
et étendue dans le cas de la dimension quelconque. Sous certaines conditions la formule
est connue (voir par exemple \cite{ks}), mais elle reste aujourd'hui à l'état de conjecture dans sa forme la plus g\'en\'erale. 
Il se trouve que la formule du conducteur de Bloch est un excellent exemple d'une question 
ouverte sur laquelle 
des progrès peuvent être obtenus à l'aide des techniques de \emph{géométrie non-commutative}
que nous présenterons par la suite. Nous pensons en réalité que cette approche non-commutative
fournira, à terme, une démonstration de la formule du conducteur de Bloch en toute généralité 
(voir 
les commentaires en fin de $\S 5$). 

Le but de cette section est de présenter cette formule, mais aussi
de rappeler ses interactions avec des formules bien connues dans le cadre de la géométrie algébrique 
(formule de type Gauss-Bonnet, de Lefschetz ou encore nombre de Milnor). \\

Une des questions typique en géométrie algébrique est l'étude de la topologie
des variétés algébriques. On se donne par exemple une variété algébrique $X$ projective sur 
$\mathbb{C}$, donnée par le lieu d'annulation d'une famille de polynômes homogènes $f_1, \dots f_p \in \mathbb{C}
[X_0,\dots,X_n]$ 
dans un espace projectif
$$X=\{x\equiv[x_0,\dots,x_n] \in \mathbb{P}^n \; | \; f_j(x)=0 \; \forall j=1, \ldots, p\}.$$
Pour simplifier on supposera que $X$ est lisse, ou en d'autre termes que le critère Jacobien 
pour $f_1,\dots,f_p$ est satisfait en tous ses points. Dans ce cas, on peut considérer 
l'espace topologique $X^{top}$ sous-jacent à $X$, dont la topologie est induite par la topologie
transcendante sur $\mathbb{C}$, qui est alors une variété topologique compacte et orientée de 
dimension $n-p$. Cette variété possède une caractéristique d'Euler bien définie
$$\chi(X):=\chi(X^{top})=\sum_{i} (-1)^idim_{\mathbb{Q}}H^i(X^{top},\mathbb{Q}) \in \mathbb{Z}.$$

Une question typique est alors de donner une formule pour calculer $\chi(X)$ en termes purement
algébrique (i.e. en ne faisant intervenir que les polynômes $f_i$ et sans faire intervenir 
la topologie transcendante sur $\mathbb{C}$). Typiquement, si $n=2$ et $p=1$, on a affaire
à une courbe algébrique plane, et on sait que l'on a
$\chi(X)=3d-d^2$
où $d$ est le degré de $f_1$.

La question possède une réponse en toute généralité de la manière suivante. Il existe sur $X$ 
une notion de $q$-formes différentielles algébriques qui forment un faisceau $\Omega_X^q$ (il
s'agit de faisceaux pour la topologie de Zariski sur $X$). Ces faisceaux possèdent des groupes de cohomologie 
$H^i(X,\Omega_X^q)$ (toujours pour la topologie de Zariski), et l'on a le théorème suivant.

\begin{thm}\label{t1}
Avec les notations ci-dessus on a 
$$\chi(X)=\sum_{i,q}(-1)^{i+q}dim_{\mathbb{C}}H^i(X,\Omega_X^q).$$
\end{thm}

Ce résultat peut se voir comme une conséquence directe de la théorie de Hodge, qui affirme en 
particulier
l'existence de décompositions
$H^k(X^{top},\mathbb{C}) \simeq \oplus_{i+q=k}H^i(X,\Omega_X^q).$
Une autre manière de comprendre cette formule est d'identifier le membre de droite 
au degré de la classe de Chern maximale $C_{top}(X)$ du fibré tangent à $X$, et d'invoquer la 
formule
de Gauss-Bonnet
$\chi(X^{top})=\int_{X}C_{top}(X).$
  
Il se trouve que la formule du Th\'eor\`eme \ref{t1} reste valable sur un corps algébriquement clos $k$ quelconque
(éventuellement de caractéristique positive), à condition d'interpréter le membre de gauche 
en termes de cohomologie $\ell$-adique. Dans ce cas, l'espace topologique $X^{top}$ ne fait plus sens
car $k$ ne porte aucune topologie naturelle intéressante, mais on peut tout de même définir 
ses groupes de cohomologie $\ell$-adique $H^*(X_{et},\Ql)$, où $\ell$ est un nombre premier différent de
la caractéristique de $k$. On définit naturellement sa caractéristique d'Euler par
$$\chi(X):=\sum_{i}(-1)^i dim_{\Ql}H^i(X_{et},\Ql).$$ 

\begin{thm}\label{t1'}
Soit $X$ une variété algébrique projective et lisse sur un corps algébriquement clos $k$. Alors on a
$$\chi(X)=\sum_{i,q}(-1)^{i+q}dim_{k}H^i(X,\Omega_X^q).$$
\end{thm}
  
Cette formule est encore une formule du type Gauss-Bonnet car le membre de droite reste 
encore égal au degré de la classe de Chern maximale $C_{top}(X) \in CH_0(X)$, considéré
comme une classe de $0$-cycles sur $X$ (ceci est une conséquence
de la formule de Hirzebruch-Riemann-Roch \cite[15.2.1]{fu} et de \cite[Ex. 3.2.5]{fu}). 
Il s'agit d'une formule typique et centrale
en géométrie algébrique, qui relie un invariant de nature topologique, $\chi(X)$,
à des invariants de nature géométrico-algébrique, à savoir les formes différentielles
algébriques. \\

Le formule du conducteur de Bloch est une formule qui s'intéresse au comportement de 
la formule du Théorème \ref{t1'} lorsque $X$ varie dans une famille 
de variétés algébriques. Pour cela, on se fixe un corps $k$, que l'on supposera parfait, 
et un anneau de valuation discrète $A$, d'uniformisante $\pi \in A$ et 
muni d'un isomorphisme $A/(\pi) \simeq k$. On supposera de plus que $A$ est un anneau hensélien 
(par exemple qu'il est complet pour la topologie $\pi$-adique). Les cas typiques à garder en tête
sont $A=k[[t]]$ avec $\pi=t$, l'anneau des séries formelles en $t$ sur $k$, ou encore
$k=\mathbb{F}_p$ et $A=\mathbb{Z}_p$ l'anneau des entiers $p$-adiques.
Le schéma $S:=Spec\, A$ 
est alors l'analogue algébrique d'un petit disque holomorphe centré au-dessus de son unique point fermé
$\{x\}=Spec\, k \subset  S$. Une famille de variétés projectives au-dessus de $S$ est par définition 
un 
schéma $X$ muni d'un morphisme projectif $X \longrightarrow S$. Concrètement, un tel schéma est 
toujours
donné par une famille de polynômes homogènes 
$F_1,\dots,F_p \in A[X_0,\dots,X_n]$, qui définissent un 
sous-schéma fermé $X \subset \mathbb{P}^n_S \longrightarrow S$ de l'espace projectif au-dessus de $S$. 

A une telle famille de variétés algébriques sur $S$ on fait correspondre deux 
variétés algébriques $X_{o}$ et $X_{t}$ définies de la manière suivante. On considère
la réduction des polynômes $F_i$ modulo $\pi$, ce qui donne des polynômes $f_i \in k[X_0,\dots,X_n]$
et donc une variété algébrique $X_{o} \subset \mathbb{P}^{n}_{\bar{k}}$ définie sur
la clôture algébrique $\bar{k}$ de $k$. De manière équivalente, on peut plonger l'anneau $A$ dans son 
corps des fractions $K=Frac(A)$, puis dans une de ses clôtures algébriques $\bar{K}$. L'image
des polynômes $F_i$ définissent ainsi des $G_i \in \bar{K}[X_0,\dots,X_p]$ et donc
une variété algébrique $X_t$ définie sur $\bar{K}$. 

Il faut rappeler ici que la situation que l'on décrit est extrêmement classique
en géométrie algébrique. Par exemple, lors de la construction des compactifications
d'espaces de modules de variétés (typiquement l'espace de modules des courbes
algébriques), on rajoute des variétés singulières comme point à l'infini, et transversalement 
autour de ce point on a exactement une famille de variétés algébriques comme-ci dessus, 
avec $X_t$ lisse et $X_o$ singulière. Un second exemple de première importance est lorsque
l'on s'intéresse aux aspects arithmétiques, par exemple avec des variétés algébriques définies
au-dessus de $\mathbb{Q}$, ou encore d'un corps de nombre. De telles variétés possèdent 
des \emph{modèles}, c'est à dire proviennent de schéma $X$ propre et
plat sur $\mathbb{Z}$ (ou sur un anneau d'entiers algébriques). 
Pour la plupart des nombre premier $p$, la réduction de $X$ modulo $p$ reste lisse, mais
il existe en général un nombre de fini de premiers pour les quels la réduction devient 
singulière. Localement autour de ce nombre premier on dispose encore d'une
famille de variétés algébriques comme ci-dessus. \\ 

Un des théorèmes fondamentaux en cohomologie étale affirme que lorsque $X_o$ et $X_t$ sont toutes 
deux lisses, 
alors on a une égalité
$\chi(X_o) = \chi(X_t).$
La formule du conducteur de Bloch se place précisément dans le cas où $X_t$ est une variété lisse sur 
$\bar{K}$ mais ou $X_o$ n'est plus nécessairement lisse sur $\bar{k}$. On s'intéresse alors à 
décrire la différence $\chi(X_o) - \chi(X_t)$, à l'aide de termes algébriques
sur le schéma $X$. Pour cela, Bloch introduit un nombre d'intersection
$deg([\Delta_X.\Delta_X]_o)$ pour le quel nous renvoyons à \cite{bl,ks} pour
une définition (qui n'est pas nécessaire pour la compréhension
de la suite de texte). Ce nombre mesure les singularités du morphisme $X \longrightarrow S$, 
et est nul quand $X_o$ est une variété non-singulière. Par ailleurs, on dispose
du conducteur de Swan $Sw(X_t)$, qui est un invariant de nature arithmétique 
de la variété $X_t$, et qui mesure les phénomènes de ramifications sauvages (ce nombre
est toujours nul en caractéristique nulle). Nous
renvoyons aussi à \cite{bl,ks} pour une définition précise, sachant que
dans ce texte nous nous placerons essentiellement 
dans le \emph{cas modéré où} $Sw(X_t)=0$. Le conjecture précise s'énonce alors comme suit.

\begin{conj}[Bloch, '85]\label{cb}
Avec les notations précédentes on a
$$\chi(X_o) - \chi(X_t) = deg([\Delta_X.\Delta_X]_o) + Sw(X_t).$$
\end{conj}

Plutôt que de tenter de détailler la nature des termes $deg([\Delta_X.\Delta_X]_o)$
et $Sw(X_t)$ dans la formule ci-dessus, nous allons rappeler plusieurs instances
de cas connus. 

\begin{enumerate}

\item Supposons que $X$ soit une variété projective et lisse sur $k$
(algébriquement clos pour simplifier). On voit
$X$ comme un schéma sur $S$ par l'application quotient $A \rightarrow k$
$$X \longrightarrow Spec\, k \longrightarrow S=Spec\, A.$$
Dans ce cas, le nombre d'intersection $deg([\Delta_X.\Delta_X]_o)$
n'est rien d'autre que l'auto-intersection de la diagonale de la variété $X$, 
qui s'écrit aussi $\sum_{i,q}(-1)^{i+q}dim_{k}H^i(X,\Omega_X^q)$. La fibre 
générique $X_t$ étant vide la formule de Bloch se réduit 
à la formule de Gauss-Bonnet discutée précédement 
$$\chi(X)=\sum_{i,q}(-1)^{i+q}dim_{k}H^i(X,\Omega_X^q).$$

\item Soit $Z$ une variété algébrique complexe lisse, et $f : Z \longrightarrow 
\mathbb{A}^1$ une fonction algébrique sur $Z$ propre. On note $A$ l'hensélisé
de l'anneau local $\OO_{\mathbb{A}^1,0}$ et $S=Spec\, A \longrightarrow \mathbb{A}^1$ le morphisme
canonique. On note $X:=Z\times_{\mathbb{A}^1}S$ le changement de base de $Z$ sur $S$. 

Supposons pour commencer que $f$ ne possède qu'une singularité isolée, ou en d'autres termes
que $X_o$ ne possède qu'un unique point singulier $x \in X_o$. On peut alors
définir l'anneau Jacobien $J(f)$ de la fonction $f$ en $x$, comme étant 
le quotient $\OO_{Z,x}/(\partial f)$, quotient de l'anneau local 
de $Z$ en $x$ par l'idéal engendré par les dérivées partielles de $f$. 
En termes plus intrinsèque, $J(f)$ est le faisceau cohérent supporté en $x$
défini comme étant le conoyau du morphisme induit par la contraction avec la différentielle
$df$
$$\xymatrix{T_X \ar[r]^-{\lrcorner \, df} & \OO_X \ar[r] & J(f).}$$
La formule du conducteur de Bloch est alors bien connue sous le nom de formule 
du nombre de Milnor et affirme que la dimension des cycles évanescents
de la singularité $x$ est égal à la dimension de l'espace $J(f)$.
Considérons la fibre
de Milnor $F_x$ de $f$ en $x$. Cette fibre de Milnor se rétracte 
par déformation sur un bouquet de sphère $\wedge_r S^{n-1}$, où $n$ est 
la dimension complexe de $X$ (voir \cite{mi}). La formule du conducteur de Bloch
s'écrit alors
$$dim_{\mathbb{C}}J(f) = r.$$

Dans le cas où les singularités de $X_o$ ne sont plus isolées, la formule reste vraie, 
et se lit de la manière suivante:
$$\chi(X_o) - \chi(X_t):=\sum_{i}(-1)^i dim_\mathbb{C}H^{i}_{DR}(X,f)$$
o\`u $H^{i}_{DR}(X,f)$ sont les groupes de cohomologie de de Rham 
de $X$ tordus par la fonction $f$ (cela est une conséquence de la suite spectrale
construite dans \cite{ka}, voir aussi \cite{sa} pour un énoncé plus précis).

\item La formule précédente à été généralisée par Deligne dans \cite[Exp. XVI]{sga7}, connue
aujourd'hui sous le nom de formule de Deligne-Milnor. Deligne démontre cette
formule dans le cas d'égale caractéristique, ce qui est équivalent à la 
formule du conducteur de Bloch lorsque l'on suppose de plus que 
$X_o$ ne possède que des singularités isolées. 

\item Enfin, la formule du conducteur de Bloch est connue lorsque l'on 
suppose que l'inclusion $X_o \hookrightarrow X$ définit le sous-schéma
réduit $(X_o)_{red}$ comme une diviseur à croisements normaux et simples
dans $X$ (voir \cite{ks}). 

\end{enumerate}

La conjecture de Bloch est aujourd'hui ouverte, y compris le cas 
minimaliste où $X_o$ est a singularités isolées en caractéristique mixte
(qui est d'une certaine façon orthogonal au cas traité dans \cite{ks}). 
Dans ce qui suit nous allons présenter une approche générale qui, nous le pensons, 
permettra d'aboutir à une preuve complète de la Conjecture \ref{cb}. Cette approche 
est basée sur des idées et techniques provenant de l'univers de la géométrie
non-commutative, et nous pensons qu'il s'agit d'un exemple pertinent 
d'interaction de la géométrie non-commutative avec des questions plus
traditionnelles de géométrie algébrique et de géométrie arithmétique. 

Avant d'entrer dans plus de détails nous souhaitons ici esquisser 
une brève idée générale de la manière dont la géométrie non-commutative
est liée à la Conjecture \ref{cb}. Il faut commencer par prendre l'exemple
$(1)$ ci-dessus au sérieux. Cet exemple est de fait relativement dégénéré, 
puisqu'il s'agit d'une famille de variétés algébriques dont la fibre générique
est vide et la fibre spéciale $X$ tout entier. C'est un cas pour lequel
la formule se réduit à la formule de Gauss-Bonnet pour $X$. Le fil directeur 
est ici de considérer ce cas comme un cas limite (une ``limite semi-classique") 
d'une situation non-commutative. En clair, nous prétendons qu'à $X \longrightarrow S$
comme dans \ref{cb} nous pouvons associer un \emph{schéma non-commutatif} $MF(X/S)$ qui 
rend compte de la théorie des singularités de $X$ au-dessus de $S$. Un point clé est que 
ce schéma non-commutatif possède une cohomologie et une caractéristique d'Euler
bien définie, et surtout, égale à $\chi(X_o)-\chi(X_t)$. Avec ce point de vue
la formule du conducteur de Bloch s'interprète essentiellement comme une
formule de Gauss-Bonnet pour le schéma non-commutatif $MF(X/S)$, qui elle même est un 
cas particulier d'une formule des traces de type Lefschetz.

Cette idée très générale peut se réaliser mathématiquement, en donnant un sens précis
à ce qu'est un schéma non-commutatif et à sa cohomologie. C'est ce que nous
allons tenter de raconter dans la suite de ce texte.

\section{\'Eléments de (non-)géométrie algébrique non-commutative}

\textbf{Un bref rappel historique.}
Il est traditionnel de considérer le théorème de Gelfand, sur la reconstruction 
d'un espace topologique de Hausdorff et compact, comme le point de départ de la géométrie non-commutative.
Pour un tel espace $X$, on note $\mathcal{C}(X)$ la $\mathbb{C}$-algèbre 
des fonctions continues sur $X$ à valeurs complexes. On considère alors 
$Spm\,  \mathcal{C}(X)$, l'ensemble des idéaux maximaux de $\mathcal{C}(X)$. Tout élément
$f \in \mathcal{C}(X)$ définit une application 
$\phi_f : Spm\,  \mathcal{C}(X) \longrightarrow \mathbb{C}$, 
qui associe à un idéal $m \subset \mathcal{C}(X)$ la classe $\overline{f} \in \mathcal{C}(X)/m \simeq
\mathbb{C}$. On munit alors l'ensemble
$Spm\,  \mathcal{C}(X)$ de la topologie la moins fine telle que 
les applications $\phi_f$ soient continues. Le théorème de reconstruction de Gelfand
s'énonce alors comme suit.

\begin{thm}[Gelfand]\label{t2}
Avec les notations précédents, l'application qui à un point $x \in X$ associe 
l'idéal $m_x \in Spm\, \mathbb{C}(X)$ des fonctions qui s'annulent en $x$, 
définit un homéomorphisme d'espaces topologiques
$$X \simeq Spm\,  \mathcal{C}(X).$$
\end{thm}

Ce résultat peut se préciser: la construction $X \mapsto \mathcal{C}(X)$ définit un 
foncteur pleinement fidèle de la catégorie des espaces separ\'es et compacts vers celle 
des $\mathbb{C}$-algèbres de Banach. Le point de départ de la géométrie non-commutative 
est alors de définir un espace non-commutatif comme étant une 
$\mathbb{C}$-algèbre de Banach, non-nécessairement commutative, qui 
ressemble à l'algèbre des fonctions sur un espace compact. Nous en renvoyons à 
\cite{co} pour plus sur ce sujet. \\

L'approche précédente à la géométrie non-commutative n'est malheureusement
plus pertinente lorsque l'on cherche à travailler dans le cadre de la géométrie
algébrique. En effet, il n'existe que très peu de fonctions sur les variétés
algébriques en général, trop peu pour espérer un énoncé analogue au Théorème
\ref{t2} lorsque $X$ est une variété algébrique. De manière plus précise, 
on voit ici deux nouveaux phénomènes concernant le comportement des fonctions
algébriques, qui ne se présentent pas dans le cadre de la topologie ou de la géométrie
différentielle.

\begin{enumerate}

\item Lorsque $X$ est une variété algébrique sur un corps algébriquement clos $k$, 
alors $X$ ne porte que des fonctions algébriques constantes dès que $X$ est 
propre (compacte) et connexe. En d'autres termes, si $\OO(X)$ désigne la $k$-algèbre 
des fonctions (algébriques) sur $X$, alors on a $\OO(X)\simeq k$. 

\item Les fonctions sur une variété algébrique $X$ s'organisent en un faisceau
d'algèbres $\OO_X$ sur $X$ (pour la topologie de Zariski). En général, ce faisceau 
possède des espaces de cohomologie $H^i(X,\OO_X)$ non-triviaux (e.g. 
lorsque $E$ est une courbe elliptique on a $H^1(X,\OO_X)\simeq k$), ce qui traduit 
des obstructions à l'extension de fonctions (manque de partition de l'unité 
par exemple). 

\end{enumerate}

Ces deux nouveaux aspects qui apparaissent dans le contexte algébrique 
obligent à adopter un point de vue différent sur la géométrie algébrique non-commutative
et à revoir l'idée, trop naïve, qu'une variété algébrique non-commutative puisse
simplement être définie comme étant une algèbre (non-nécessairement commutative).

Concernant le point $(1)$ ci-dessus, il faut commencer par remarquer que bien qu'il 
n'existe que très peu de fonctions définies globalement, il existe en revanche de nombreuses
sections de fibrés vectoriels algébriques. Cela suggère d'élargir l'algèbre
des fonctions $\OO(X)$ en considérant tous les fibrés vectoriels sur $X$ ainsi que
tous les morphismes entre de tels fibrés. Ces fibrés forment une catégorie
$Vect(X)$, qui est $k$-linéaire si de plus $X$ est définie sur le corps $k$. Par ailleurs,
$Vect(X)$ contient l'algèbre $\OO(X)$ comme endomorphisme du fibré trivial 
de rang $1$ sur $X$. Les aspects cohomologiques du point $(2)$ peuvent quand 
à eux être considérés et faisant de $Vect(X)$ une \emph{dg-catégorie} $k$-linéaire. 
Il s'agit d'une structure de type catégorie où les morphismes entre deux
objets forment des complexes de $k$-modules. Pour deux fibrés vectoriels $E$ et $F$ sur
$X$ il est possible de construire un complexe naturel $\mathbb{R}\underline{Hom}(E,F)$
dont la cohomologie calcule les $Ext$ globaux entre $E$ et $F$ et qui organisent 
$Vect(X)$ en une dg-catégorie. En conclusion, ce que nous dicte cette discussion est 
qu'une variété algébrique non-commutative peut, et doit, être définie comme
une dg-catégorie $k$-linéaire. C'est un point de vue largement adopté aujourd'hui, 
par exemple par des auteurs tels que Bondal, Kontsevich Kapranov, Orlov, Van den Bergh ...
  \\

\textbf{La notion de schémas non-commutatifs.}
Entrons un peu plus dans les détails. On se fixe un anneau commutatif de base $k$ (qui pourra
par la suite être notre corps de base, mais il est important de s'autoriser
ce degré de généralité pour la suite). Par définition, une dg-catégorie $T$ sur $k$
consiste en les données suivantes. 

\begin{enumerate}

\item Un ensemble $Ob(T)$, que l'on appelle l'ensemble des objets de $T$.

\item Pour toute paire d'objets $(x,y) \in Ob(T)\times Ob(T)$
un complexe de $k$-modules $T(x,y)$, que l'on appelle le complexe
des morphismes de $x$ vers $y$.

\item Pour tout triplet d'objets $(x,y,z) \in Ob(T)\times Ob(T)\times Ob(T)$, 
un morphisme de complexes
$$T(x,y) \otimes_k T(y,z) \longrightarrow T(x,z),$$
appelé composition. 

\end{enumerate}

On demande par ailleurs que les morphismes de compositions satisfassent à des 
conditions d'associativité et d'unité naturelles (voir \cite{crm} pour les détails). En utilisant un langage 
différent, une dg-catégorie est une catégorie enrichie dans la catégorie monoidale
$(C(k),\otimes_k)$ des complexes de $k$-modules. 

\begin{df}\label{d1}
Une \emph{schéma non-commutatif}, ou encore \emph{nc-schéma}, (au dessus de $k$) est une 
dg-catégorie sur $k$. 
\end{df}

La définition ci-dessus peut sembler naïve à priori, et pour tout dire 
il tient du miracle qu'une notion aussi générale et simple se révèle à postériori pertinente. 
Le lecteur intéressé pourra consulter \cite{crm} qui tente d'expliquer pourquoi 
les dg-catégories sont plus pertinentes que les simples catégories (triangulées par exemple). 

Pour appréhender la notion de nc-schémas ci-dessus il faut commencer par expliquer
comment elle est reliée à la notion de variété algébrique ou plus généralement de schéma. 
Soit donc une variété algébrique $X$ sur un corps $k$, 
que l'on va supposer être affine pour commencer. 
Elle dispose d'une $k$-algèbre de fonctions $A:=\OO(X)$. On associe à cette 
algèbre $A$ une dg-catégorie $\D(A)$ sur $k$ de la manière suivante. Les objets
de $\D(A)$ sont par définition les complexes bornés de 
$A$-modules projectifs de type fini. Pour deux tels complexes $E$ et $F$, on 
dispose d'un complexe de morphismes $A$-linéaires $\underline{Hom}_A(E,F)$, qui 
est naturellement un complexe de $k$-module. Cela définit les 
complexes de morphismes dans $\D(A)$, qui munit de la composition 
usuelle des morphismes définit la dg-catégorie $\D(A)$. 

Supposons maintenant que $X$ ne soit plus nécessairement affine. Dans ce cas, on recouvre
$X$ par un nombre fini de sous-variétés ouvertes affines $U_i \subset X$ pour les quelles
on sait définir $\D(U_i)$. On construit alors $\D(X)$ par un procédé de recollement 
à l'aide des $\D(U_i)$ (voir \cite{mf}). 
Il faut mettre en garde ici sur le fait que ce procédé de recollement
est relativement subtil et fait intervenir des outils d'algèbre homotopique (voir par exemple
\cite[5.3.1]{crm}). 
La construction précédente n'utilise nullement le fait que $X$ soit une variété 
algébrique et reste valable lorsque $k$ est un anneau commutatif et $X$ est un $k$-schéma.
De manière plus concrète $\D(X)$ peut aussi se définir comme la dg-catégorie
des complexes de $\OO_X$-modules injectifs, bornés à gauche, et localement quasi-isomorphes
à des complexes bornés de $\OO_X$-modules libres de rang fini (voir \cite[\S 5.3]{crm}
où cette dg-catégorie est notée $L_{pe}(X)$). Mais plutôt que d'entrer
dans les détails définitionnels de $\D(X)$ citons quelques unes de ses propriétés clés. 

\begin{itemize}

\item Si $V$ est un fibré vectoriel sur $X$, que l'on voit comme un faisceau de 
$\OO_X$-modules localement libres, alors $V$ définit un objet de $\D(X)$ en considérant
$V$ comme un complexe de manière triviale (valant $V$ en degré nul, et $0$ en tout autre
degré). Plus généralement on peut décaler $V$ pour le placer en un unique degré $n$
pour obtenir un nouvel objet $V[-n] \in \D(X)$. La construction $V \mapsto V[0]$, 
induit un foncteur (pleinement fidèle en un certain sens) 
de la catégorie des fibrés vectoriels sur $X$ vers
$\D(X)$. 

\item Pour deux fibrés vectoriels $V$ et $W$, et tout entier $n$ 
on dispose d'isomorphismes naturels
$$H^i(\D(X)(V,W[n])) \simeq Ext^{i+n}_{\OO_X}(V,W).$$
En particulier, on a $H^i(\D(X)(\OO_X,W)) \simeq H^i(X,W)$. 

\end{itemize}

Les deux propriétés ci-dessus peuvent être comprises comme caractérisant $\D(X)$ en un 
certain sens, et dans la pratique elles sont les deux outils essentiels pour
appréhender $\D(X)$. On voit bien par ailleurs qu'elles répondent précisément aux
problématiques que nous avions soulevées plus haut. \\

La construction $X \mapsto \D(X)$ est une première source d'exemples de nc-schémas. Il existe
beaucoup d'autres, d'origines variées. 

\begin{enumerate}

\item Pour commencer, on trouve de nombreux exemples intéressants de nc-schémas
comme "morceaux" de variétés algébriques. Pour une variété $X$, il arrive
que sa catégorie dérivée $\D(X)$ se décompose comme une produit semi-direct
de deux sous-catégories $\mathcal{A}, \mathcal{B} \subset \D(X)$ (appelé
\emph{décomposition semi-orthogonale}). Les morceaux $\mathcal{A}$ et 
$\mathcal{B}$ sont des nc-schémas qui ne sont généralement plus 
de la forme $\D(Y)$ pour une variété $Y$. L'étude des décompositions
semi-orthogonales des variétés projectives lisses est par exemple intimement 
liée à des questions profondes de rationalité. Pour être plus précis, 
on peut montrer qu'une cubique $X$ de dimension $4$ possède un facteur
$\mathcal{A} \subset \D(X)$ qui est une surface K3 non-commutative 
Par ailleurs, Kuznetsov conjecture que $X$ est rationnelle si et seulement si 
$\mathcal{A}$ est une K3 commutative (i.e. de la forme $\D(S)$ pour une
surface K3 $S$, voir \cite[Conj. 1.1]{ku2}). L'étude des facteurs directs de $\D(X)$ du point de vue
de la géométrie non-commutative est un sujet en soi et aujourd'hui très actif
(voir \cite{ku1}).

\item Si $A$ est une $k$-algèbre, associative et unitaire
mais non nécessairement commutative, on dispose aussi d'une dg-catégorie $\D(A)$ 
formée des complexes bornés de $A$-modules projectifs de type fini. Il faut penser à
$\D(A)$ comme au schéma non-commutatif ``$Spec\, A$". Un  exemple important est 
lorsque $A=k[\Gamma]$ est la $k$-algèbre en groupes d'un groupe discret $\Gamma$. 

\item Si un groupe algébrique $G$ opère sur une variété alégbrique $X$, on dispose
d'un nc-schéma quotient de $X$ par $G$ que l'on note $\D^G(X)$. La dg-catégorie 
$\D^G(X)$ est définie en considérant les objets $G$-équivariants de $\D(X)$ et le nc-schéma
correspondant doit être vu comme le quotient non-commutatif de $X$ par l'action de $G$. 

\item Certains nc-schémas arrivent comme déformations de schémas. Un exemple typique
provient de la quantification par déformation des variétés algébriques symplectiques (ou de
Poisson). Pour une telle variété $X$ la déformation par quantification s'incarne en 
une dg-catégorie $k[[t]]$-linéaire $\D^q(X)$, qui est une déformation formelle
de $\D(X)$. 

\end{enumerate}

\medskip

Pour clore cette discussion sur les nc-schémas nous signalons qu'il existe, en quelque sorte, 
un analogue au théorème de Gelfand. Pour une variété algébrique $X$ sur $k$, ou plus
généralement un $k$-schéma que l'on supposera quasi-compact et quasi-séparé, 
on peut montrer que la dg-catégorie $\D(X)$ est engendrée par un unique objet $E$ (voir 
\cite{bv}). 
Un tel objet est appelé un générateur compact, et son existence, qui est un énoncé
d'existence relativement profond, possède de nombreuses implications. Tout d'abord, 
on peut considérer les endomorphismes de l'objet $E$, qui munient de la composition
des endomorphismes forment une dg-algèbre $B_X:=\D(X)(E,E)$. Le fait que $E$ soit un générateur
compact de $\D(X)$ s'interpréte aussi par l'existence d'une équivalence de dg-catégories
$$\D(X) \simeq \D(B_X),$$
où $\D(B_X)$ est la dg-catégorie des $B_X$-dg-modules parfaits. La dg-algèbre 
$B_X$ est en quelque sorte l'analogue de l'algèbre des fonctions $\mathcal{C}(X)$
du théorème de Gelfand \ref{t2}. Cependant, la situation est ici un peu différente. 

\begin{enumerate}

\item En général, la dg-algèbre $B_X$ ne permet pas de reconstruire la variété (ou le schéma)
$X$. En effet, il existe de nombreux exemples de variétés $X$ et $Y$ non isomorphes, mais
telles ques $\D(X) \simeq \D(Y)$. Dit autrement, deux variétés algébriques non isomorphes peuvent
être isomorphes en tant que schémas non-commutatifs ! Il y a des raisons de penser que
la relation d'équivalence que cela introduit sur les variétés algébriques est une relation
pertinente, par exemple pour les questions de géométrie birationnelle (voir \cite{ku2}). 

\item La dg-algèbre $B_X$ n'est pas canonique et dépend du choix d'un générateur
compact $E \in \D(X)$. En particulier, $B_X$ n'est pas fonctoriel en $X$, seule la dg-catégories 
$\D(B_X)$ l'est. En d'autres termes, il faut considérer $B_X$ uniquement à équivalence
de Morita près (voir \cite[\S 4.4, Def. 8]{crm}).

\end{enumerate}

\textbf{La catégorie des nc-schémas.} A ce stade il nous faut 
dire un mot ou deux sur la notion de morphismes entre nc-schémas sur la quelle nous sommes
restés silencieux pour l'instant. Malheureusement c'est précisément le point où les choses
se compliquent considérablement et pour le quel il est nécessaire d'utiliser 
le langage des infinies-catégories (voir \cite[\S 2.1]{ems}). Nous ne pouvons raisonnablement pas 
faire un détour vers les $\s$-catégories dans ce texte, et nous nous contenterons
donc de renvoyer vers d'autres références. Il est possible de remplacer 
formellement l'expression "$\s$-catégorie" par "catégorie" dans ce qui suit,
cependant il faut garder en tête qu'alors plusieurs énoncés seront tout simplement incorrects. 

Il existe tout d'abord une notion évidente de morphisme entre dg-catégories que nous appellerons
morphisme strict. Un tel morphisme $f : T \longrightarrow T'$ entre deux dg-catégories
$T$ et $T'$ est la donnée d'une application $Ob(T) \rightarrow Ob(T')$ entre les ensembles
d'objets, et pour toute paire $(x,y)$ d'objets de $T$ d'un morphisme de complexes
$k$-linéaires $f_{x,y} : T(x,y) \longrightarrow T'(f(x),f(y))$. On demande par ailleurs que
les morphismes $f_{x,y}$ respectent les compositions et les unités de $T$ et $T'$. 
Comme son nom l'indique, la notion de morphisme strict n'est malheureusement pas adaptée
car trop rigide. On souhaite en effet identifier certaines dg-catégories qui ne sont
pas strictement isomorphes. Pour cela il faut introduire la notion d'\emph{équivalence Morita}.
Il s'agit, en gros des morphismes stricts 
$f : T \longrightarrow T'$, qui vérifient les deux conditions suivantes.
\begin{itemize}

\item Tous les morphismes de complexes $f_{x,y} : T(x,y) \longrightarrow T'(f(x),f(y))$
sont des quasi-isomorphismes (i.e. induisent des isomorphismes sur les groupes
de cohomologie correspondants). 

\item L'image de $T$ par $f$ \emph{engendre} $T'$ par les opérations
qui consistent à prendre un facteur direct, le cône
d'un morphisme ou le décalage.

\end{itemize}

La notion d'équivalence Morita entre dg-catégorie est la notion pertinente dans notre contexte, 
et il faut donc forcer le fait que ces équivalences sont des morphismes inversibles. Cela se fait par
le procédé de localisation des catégories: on considère la catégorie des dg-catégories
et morphismes stricts à laquelle on rajoute de manière formelle des inverses à toutes les
équivalences Morita (voir \cite[\S 4.4]{crm}). Il se trouve que ce procédé de localisation est relativement profond, 
et le résultat n'est en général pas une catégorie mais une \emph{$\s$-catégorie}, comme cela est par
exemple expliqué dans \cite[\S 2.1]{ems}. 

L'$\s$-catégorie ainsi obtenue, en localisant les dg-catégories le long des équivalence Morita, 
est, par définition, l'opposée de l'$\s$-catégorie des schémas non-commutatifs sur $k$. On notera
$$k-\Sch^{nc}:=(dg-cat_k)[Morita^{-1}]^{op}.$$
Les objets de l'$\s$-catégorie $k-\Sch^{nc}$ sont simplement les dg-catégories sur $k$. En revanche les
morphismes ne sont pas faciles à décrire. Un morphisme strict induit évidemment un morphisme
dans $k-\Sch^{nc}$. Réciproquement on peut montrer que tout morphisme $T \longrightarrow T'$  dans 
$k-\Sch^{nc}$ peut se représenter par un diagramme de morphismes stricts
$\xymatrix{
T & \ar[l]_-{v} \ar[r]^-{u} T'' & T',}$
avec $v$ une équivalence Morita. La combinatoire possible de ce type de diagrammes
s'exprime précisément dans la structure $\s$-catégorique portée par $k-\Sch^{nc}$. \\

\textbf{Une géométrie des nc-schémas ?} Nous venons de voir l'existence d'une $\s$-catégorie
des nc-schémas. Il est très certainement naturel de se poser la question de l'existence
d'une \emph{géométrie des nc-schémas}. Cette question est évidemment floue, mais est 
souvent interprétée comme la question de savoir étendre aux nc-schémas des notions
standards de géométrie des variétés algébriques et des schémas. Nous souhaitons signaler
ici qu'il n'existe pas, à proprement parler, de géométrie des nc-schémas. 
On peut s'en convaincre en considérant deux notions fondamentales en géométrie
des schémas: la notion d'ouvert de Zariski et celle de points. Comme nous allons le voir, 
aucune de ces deux notions ne possède d'extension pertinente au cadre des nc-schémas.

\begin{itemize}

\item Les points non-commutatifs. Pour une variété algébrique $X$ sur un corps $k$
(ou plus généralement pour un $k$-schéma), les points (rationnels sur $k$) de $X$
sont exactement les morphismes $Spec\, k \longrightarrow X$. Si l'on adopte 
cette définition pour les nc-schémas on trouve qu'un point d'un nc-schéma correspondant à
une dg-catégorie $T$ est un objet de $T$ (on suppose ici que $T$ est triangulée au sens 
de \cite[\S 4.4]{crm}). 
L'ensemble des objets d'une dg-catégorie $T$, même si l'on impose de fortes
conditions de finitude sur $T$, tendent à former des espaces relativement pathologiques. 
Ces espaces sont en général très fortement non-séparés et présentent un nombre dénombrable
de composantes irréductibles, y compris lorsque l'on suppose $T$ de la forme $\D(X)$
pour $X$ une variété propre et lisse (voir \cite{tv}). Il parait délicat, voire impossible, 
de penser faire de la géométrie (au sens l'on en fait avec des variétés algébriques)
avec cette notion de points.

\item Les ouverts non-commutatifs. La situation avec les ouverts de Zariski est 
relativement similaire. Pour un schéma $X$, tout ouvert $U \subset X$ peut-être
vu comme le complémentaire du support d'un complexe parfait sur $X$. Ce fait 
incite naturellement à définir les ouverts de nc-schémas comme étant 
les complémentaires d'objets. Plus précisément, si un nc-schéma est donnée par
une dg-catégorie $T$, et si $K \in T$ est un de ses objets, on peut 
considérer le morphisme quotient $T \longrightarrow T/<K>$, qui consiste
à rendre nul l'objet $K$ de manière universelle. Ce morphisme quotient 
définit un morphisme de nc-schémas $T/<K> \longrightarrow T$ 
que l'on voit comme un ouvert de Zariski de $T$ (complémentaire du support de $K$).
Cependant, si l'on adopte cette définition, on voit que la droite projective 
$\mathbb{P}^{1}$ possède un ouvert isomorphe à un point $Spec\, k$ (e.g. comme
complémentaire de l'objet $\OO_{\mathbb{P}^1}$).

\end{itemize}

\medskip

Les deux exemples précédents montrent que les extensions naturelles des notions de points
et d'ouverts au cadre des nc-schémas ne sont pas pertinentes, tout au moins si l'on souhaite
faire de la géométrie en un sens relativement standard. Il n'y a, à priori, pas de raisons
formelles à ce que d'autres définitions, plus pertinentes, puissent exister. Cependant, nous
sommes convaincus que chercher à faire de la géométrie avec des nc-schémas, tout en cherchant à 
rester proche de la géométrie des schémas et des variétés, 
est probablement voué à l'échec. De manière plus brutale et polémique, 
nous pensons que la géométrie des nc-schémas n'existe pas ! Ce point de vue est par ailleurs 
partagé, et a par exemple amené Kontsevich à introduire l'expression amusante 
de \emph{non-commutative
non-geometry}. 

D'une certaine manière, c'est la nature de la théorie des nc-schémas de ne pas préserver les
notions standards de la géométrie algébrique, tels que les ouverts, les points, 
ou encore des notions plus avancées comme la platitude. La philosophie sous-jacente est ici 
qu'il faut accepter de perdre sur certains aspects pour gagner sur d'autres. Ici, on gagne
très certainement sur la souplesse et la simplicité de la définition de nc-schémas, qui
permet de son côté de produire de très nombreux exemples. D'un autre côté, nous allons voir 
que les nc-schémas possèdent des théories cohomologiques forts pertinentes, y compris
des notions sophistiquées telles la cohomologie $\ell$-adique. Ce fait est des plus surprenant 
sachant que la notion de topologie (que se soit de Zariski ou étale) n'est plus
disponible pour les nc-schémas. On touche là la force de la théorie des nc-schémas: 
la définition de nc-schéma est extrêmement générale (presque naïve à première vue), mais
on sait tout de même définir la cohomologie de tels objets, et donc certains invariants
numériques de type caractéristique d'Euler. 

\section{Cohomologie des nc-schémas}

Nous venons de voir l'existence d'une $\s$-catégorie de schémas non-commutatifs $k-\Sch^{nc}$, 
munie d'un foncteur $k-\Sch \longrightarrow k-\Sch^{nc}$ qui associe à tout $k$-schéma $X$ 
le schéma non-commutatif correspondant $\D(X)$. Nous avons aussi vu que les notions usuelles
de géométrie ne sont plus toujours pertinentes dans le monde non-commutatif. Cependant, certaines
constructions et outils de la géométrie algébrique se généralisent au cadre non-commutatif, 
et c'est le cas en particulier des invariants de nature cohomologique. \'Etant donné
le manque de notion de géométrie de base dans le cadre non-commutative (essentiellement
absence de la notion de topologie) ceci peut paraitre surprenant (et ça l'est aux yeux
des auteurs). Il faut cependant relativiser: tout se que l'on sait faire
avec la cohomologie des schémas ne s'étend pas au cas des nc-schémas. L'essentiel pour nous 
est que la notion de caractéristique d'Euler garde un sens, ainsi que la formule de traces
qui permet de la décrire en termes d'invariants algébriques. C'est ce que nous allons
décrire maintenant. 

Pour une variété algébrique propre et lisse $X$ sur un corps $k$, on dispose 
des espaces de cohomologie $\ell$-adique $H^i_{et}(\overline{X},\Ql)$. Ce sont
des $\Ql$-espaces vectoriels de dimension fini sur le quel le groupe
de Galois $Gal(k^{sp}/k)$ opère. Par ailleurs, on dispose aussi
des faisceaux des formes différentielles $\Omega_X^q$ et de leur 
cohomologie $H^i(X,\Omega_X^q)$. Ces groupes de cohomologie ne se comportent pas
vraiment de manière raisonnable et ne forment pas une bonne théorie cohomologique
à proprement parler. Cependant, on dispose de la formule de Gauss-Bonnet (théorème \ref{t1'})
 qui 
prédit que les deux caractéristiques d'Euler correspondantes sont égales. 
Il se trouve que cette situation persiste dans le contexte non-commutatif. Dans cette section
nous présentons les versions non-commutatives de ces deux types de cohomologie. \\

\textbf{Homologie de Hochschild comme formes différentielles non-commutatives.}
Commençons par le cas de la cohomologie de Hodge $H^i(X,\Omega_X^q)$. La version non-commutative
de cette cohomologie est connue depuis des décennies et s'appelle l'homologie de Hochschild. Historiquement elle a d'abord
été introduite par Hochschild 
dans le cadre des algèbres puis généralisée à des cadres plus généraux dont 
celui des dg-catégories (voir par exemple \cite[\S 5.3]{ke}). 
Le point important pour nous est que l'homologie de Hochschild
possède en réalité une interprétation en termes de traces dans une $\s$-catégorie monoidales
adéquate, ce que nous verrons plus tard dans la section suivante. Nous nous contentons ici de 
rappeler
brièvement sa construction. 

Soit donc un nc-schéma sur un anneau $k$, donné par une dg-catégorie $T$. On construit un complexe
de Hochschild $C_*(T)$. Ce complexe est un modèle à ce que l'on est en droit de noter
$T\otimes^{\mathbb{L}}_{T\otimes^{\mathbb{L}}_k T^o}T$. Les formules explicites
décrivant le complexe $C_*(T)$ se trouvent par exemple dans \cite[\S 5.3]{ke}.
Le cohomologie du complexe $C_*(T)$ est par définition l'homologie de Hochschild de $T$, que 
l'on peut
aussi appeler la cohomologie de Hodge du schéma non-commutatif correspondant. On note ces 
groupes
par $HH_i(T):=H^{-i}(C_*(T))$, qui sont naturellement des $k$-modules. En général, sans 
hypothèses additionnelles
sur $T$ ces groupes sont des $k$-modules potentiellement arbitraires et ne satisfont aucune
condition de finitude. Nous verrons à la section suivante qu'une condition simple sur $T$
(propreté et lissité au sens non-commutatif) assure que $HH_i(T)$ sont de type fini, et
plus généralement que la caractéristique d'Euler $\sum (-1)^i[HH_i(T)]$ est bien définie
comme élément de $K_0(k)$. Pour l'instant nous nous contenterons des quelques propriétés 
importantes suivantes.

\begin{enumerate}

\item Soit $X$ une variété algébrique lisse sur un corps $k$, et considérons 
le nc-schéma correspondant $\D(X)$. On dispose alors d'isomorphismes naturels 
(appelé HKR pour Hochschild-Kostant-Rozenberg)
de $k$-espaces vectoriels (pour $i \in \mathbb{Z}$)
$$HH_i(\D(X)) \simeq \bigoplus_{p-q=i}H^p(X,\Omega^q_X).$$

\item Si $A$ est une $k$-algèbre plate, et $\D(A)$ le nc-schéma correspondant, alors
on a 
$$HH_i(\D(A)) \simeq Tor_i^{A\otimes_k A^o}(A,A).$$

\item La construction $T \mapsto C_*(T)$ envoie équivalences Morita de dg-catégories
sur quasi-isomorphismes de complexes. Elle induit en particulier un $\s$-foncteur
$C_* : k-\Sch^{nc} \longrightarrow L(k)^{op},$
où $L(k)$ est l'$\s$-catégorie des complexes de $k$-modules (obtenue en localisant 
la catégorie des complexes le long des quasi-isomorphismes, voir \cite{crm}).

\end{enumerate}

Il existe aussi une version à coefficients. Pour $T$ une dg-catégorie et 
$f : T \longrightarrow T$ un endomorphisme de $T$, on dispose d'un complexe
$C_*(T,f)$, qui est cette fois un modèle à $T\otimes^{\mathbb{L}}_{T\otimes^{\mathbb{L}}_k T^o}
\Gamma(f)$. Ici $\Gamma(f)$ est le graphe de $f$, c'est à dire le bimodule sur $T$
qui envoie $(x,y)$ sur $T(y,f(x))$. Lorsque $f=id$, on retrouve le complexe de Hochschild ci-
dessus $C_*(T,id)=C_*(T)$. Il faut penser que $C_*(T,f)$ représente 
l'intersection entre la diagonale de $T$ et le graphe de $f$, et est donc 
un avatar algébrique des points fixes de l'endomorphisme $f$. En particulier,
lorsqu'elle est définie, la caractéristique d'Euler de $C_*(T,f)$ 
est en droit de s'appeler le \emph{nombre de Lefschetz de $f$}, qui représente
le nombre virtuel de points fixes. \\

\textbf{Cohomologie $\ell$-adiques des nc-schémas.} Nous venons de voir que l'homologie de 
Hochschild  des dg-catégories permettait de donner un sens à la cohomologie de Hodge des 
nc-schémas. Le cas de la cohomologie $\ell$-adique est une autre histoire et sa définition dans 
le contexte non-commutatif est très récente. Elle est basée sur des résultats de Thomason 
qui affirment que la cohomologie $\ell$-adique d'un schéma peut être reconstruite à partir
de sa K-théorie algébrique. Cette idée à été reprise par Blanc dans \cite{bla} 
afin d'introduire
la notion de cohomologie de Betti rationnelle de dg-catégories complexes. Dans \cite{mf}, 
cette 
construction est reprise dans le contexte $\ell$-adique. Les détails de la construction
dépassent le cadre de ce texte de survol, nous allons nous contenter d'esquisser l'idée 
générale et de rappeler les principales propriétés. 

Soit donc $T$ une dg-catégorie sur $k$ représentant un nc-schéma. Pour simplifier, mais ce 
n'est pas strictement nécessaire, nous supposerons que $T$ ne possède qu'un unique objet, et 
est donc donnée par une dg-algèbre $B$. L'idée est alors d'approximer
$T$ par des $k$-schémas commutatifs affines lisses sur $k$ de la manière suivante. Pour tout
schéma affine lisse $Spec\, A$ sur $k$, on considère la dg-algèbre $A \otimes_k B$. 
On regarde $K(A\otimes_k B)$, l'espace de K-théorie algébrique de la dg-algèbre $A\otimes_k B$,
L'association $A \mapsto K(A\otimes_k B)$ défini un foncteur 
des schémas affines lisses sur $k$ vers celle des espaces. Ce foncteur sera noté
$K^T$ et définit une théorie cohomologique généralisée pour les $k$-schémas. En d'autres
termes $K^T$ définit un objet de  $\SH_k$, la ($\s$-)catégorie homotopique stable des $k$-schémas
au sens de Morel-Voevodsky (voir \cite{mv}). 
On observe que lorsque $T=k$ est la dg-catégorie unitaire (un seul objet et $k$ comme
endomorphisme), alors $K^k=\BU_k$ est le spectre de K-théorie motivique au-dessus de $k$
(aussi noté $KGL$ dans la littérature). Cela permet de voir que 
$K^T$ est naturellement muni d'une structure de modules au-dessus de $\BU_k$. On utilise alors
la construction de réalisation $\ell$-adique de \cite{ay}, 
qui est par définition un $\s$-foncteur
$r_{\ell} : \SH_k \longrightarrow \D(S,\Ql),$
où $\D(S,\Ql)$ est l'$\s$-catégorie des complexes $\Ql$-adiques ind-constructibles sur 
$S:=Spec\, k$.  La réalisation de $K^T$ est ainsi un objet de $\D(S,\Ql)$, qui est par ailleurs
un module sur $r_\ell(\BU_k)$. Il est facile de voir que l'on a
$$r_{\ell}(\BU_k)\simeq \Ql(\beta):=\bigoplus_{i \in \mathbb{Z}}\Ql[2i](i).$$

\begin{df}\label{d2}
La \emph{cohomologie $\ell$-adique de $T$} est 
le $\Ql(\beta)$-module $r_{\ell}(K^T)$ défini ci-dessus. Il sera noté
$$\mathbb{H}(T,\Ql):=r_{\ell}(K^T) \in \Ql(\beta)-Mod.$$
\end{df}

Notons que par définition $\mathbb{H}(T,\Ql)$ est un complexes de faisceaux de $\Ql$-espaces
vectoriels sur $S_{et}$, le petit site étale de $S=Spec\, k$. Notons aussi que
cet objet est un module sur $\Ql(\beta)$, ou en d'autres termes il est \emph{$2$-périodique à 
twist de Tate près}
$$\mathbb{H}(T,\Ql)[2](1) \simeq \mathbb{H}(T,\Ql).$$
Donc, lorsque $k$ contient les racines de l'unité, ceci se transforme en un 
$2$-périodicité $\mathbb{H}(T,\Ql)[2] \simeq \mathbb{H}(T,\Ql)$. \\

Les deux propriétés fondamentales de la 
construction $T \mapsto \mathbb{H}(T,\Ql)$ sont les suivantes.

\begin{enumerate}

\item Supposons que $p : X \longrightarrow S=Spec\, k$ soit un $k$-schéma séparé et de type 
fini. On supposera de plus
que soit le morphisme $p$ est propre, soit $k$ est un corps. Alors, 
on dispose d'isomorphismes naturels dans $\D(S,\Ql)$
$$\mathbb{H}(\D(X),\Ql) \simeq \mathbb{R}p_*(\Ql(\beta)) \simeq 
\bigoplus_{i \in \mathbb{Z}}\mathbb{R}p_*(\Ql)[2i](i).$$
En particulier, si $k$ est un corps algébriquement clos on a un
isomorphisme de $\Ql$-espaces vectoriels 
$$\mathbb{H}^k(\D(X),\Ql) \simeq 
\bigoplus_{i \in \mathbb{Z}}H_{et}^{k+2i}(X,\Ql)$$

\item Si $T_0 \subset T$ est une sous-dg-catégorie, et si 
$T/T_0$ est le quotient dans l'$\s$-catégorie des dg-catégories 
et morphismes Morita, alors on a un triangle distingué
$$\xymatrix{
\mathbb{H}(T_0,\Ql) \ar[r] & \mathbb{H}(T,\Ql) \ar[r] & \mathbb{H}(T/T_0,\Ql) \ar[r]^-{+1} 
& }$$
En particulier, on dispose d'une suite exacte longue en cohomologie
$$\xymatrix{
... \ar[r] & \mathbb{H}^k(T_0,\Ql) \ar[r] & \mathbb{H}^k(T,\Ql) \ar[r] & \mathbb{H}^k(T/T_0,
\Ql) \ar[r] & 
\mathbb{H}^{k+1}(T_0,\Ql) \ar[r] & ...}$$

\end{enumerate}

\section{Une formule des traces}

Nous venons de voir qu'il était possible de définir la cohomologie des schémas non-commutatifs, 
et de fait que du point de vu cohomologique les nc-schémas partageaient des similarités avec 
les schémas. Nous allons maintenant observer un phénomène purement non-commutatif, qui n'a pas
d'analogue pour les schémas: la dualité. Appliquer à un nc-schéma de la forme $\D(X)$ pour
$X$ une variété propre et lisse, cette dualité est une incarnation de la dualité de Poincaré. 
Un point clé est qu'ici cette dualité existe déjà dans la catégorie des schémas 
non-commutatifs, avant d'en prendre la cohomologie. 

On commence par observer que l'$\s$-catégorie $k-\Sch^{nc}$ des nc-schémas
est munie d'une structure monoïdale symétrique $\otimes_k$, qui est la version
non-commutative du produit de deux schémas. Pour deux nc-schémas $T$ et $T'$, leur
produit tensoriel $T\otimes_k T'$ est simplement la dg-catégorie produit tensoriel
de $T$ avec $T'$ au-dessus de $k$. L'ensemble des objets de $T\otimes_k T'$ est 
l'ensemble produit $Ob(T) \times Ob(T')$, et pour deux couples $(x,x')$ et $(y,y')$ 
le complexe des morphismes est le produit tensoriel des complexes de morphismes
$T(x,y) \otimes_k T'(x',y')$. Cette opération munit $k-\Sch^{nc}$ d'une structure monoïdale
symétrique $\otimes_k$, qui se trouve être compatible avec la structure produit sur la
catégorie des $k$-schémas: l'$\s$-foncteur $X \mapsto \D(X)$ se promeut en un $\s$-foncteur
monoïdal symétrique (c'est une conséquence des résultats de \cite{to3})
$\D : k-\Sch \longrightarrow k-\Sch^{nc}.$

Rappelons que pour une catégorie monoïdale symétrique $(C,\otimes)$, on dispose d'une notion 
d'objet \emph{dualisable}, qui formalise la notion de dualité parfaite
entre espaces vectoriels de dimension finie par exemple. Un objet $x \in C$ est dit dualisable
s'il existe un objet $x^\vee \in C$ et deux morphismes (appelés coévaluation et évaluation)
$$\xymatrix{
\mathbf{1} \ar[r]^-{coev} & x \otimes x^{\vee} \ar[r]^-{ev} & \mathbf{1},}$$
qui vérifient les identités triangulaires naturelles. On montre que lorsque
un dual $x^\vee$ existe il est automatiquement unique à isomorphisme unique près. \^Etre 
dualisable est ainsi une propriété et les données $(x^\vee,ev,coev)$ sont entièrement 
déterminées par l'objet $x$ lorsqu'elles existent. Par exemple, lorsque
$(C,\otimes)$ est la catégorie des modules sur une anneau commutatif $k$, on montre facilement
que les objets dualisables sont exactement les $k$-modules projectifs et de type fini. 
Pour un tel module $M$ son dual $M^{\vee}$ est canoniquement isomorphe au module
$\underline{Hom}(M,k)$ et l'application $ev$ est l'application d'évaluation 
naturelle qui envoi $m\otimes u$ sur $u(m)$. 

Ce que nous venons de résumer pour les catégories garde un sens pour les $\s$-catégories
(voir par exemple \cite{chern}), et nous posons alors la définition suivante.

\begin{df}\label{d3}
Un nc-schéma $T \in k-\Sch^{nc}$ est \emph{saturé} s'il est dualisable 
dans l'$\s$-catégorie monoïdale symétrique $(k-\Sch^{nc},\otimes_k)$.
\end{df}

Il se trouve que l'on peut décrire explicitement ce que sont les nc-schémas saturés. On montre 
en effet (voir \cite[Prop. 2.5]{to4}) 
qu'il s'agit exactement des dg-algèbres \emph{propres et lisses}. Rappelons
qu'une dg-algèbre $B$ sur $k$ est propre si le complexe sous-jacent à $B$ est parfait (i.e.
quasi-isomorphe à un complexe borné de $k$-modules projectifs de type fini). Elle est dite
lisse si l'objet $B$ possède une résolution finie en tant que $B\otimes_k B^o$-dg-module
(de manière plus précise si c'est un objet compact de la catégorie
dérivée des bimodules $D(B\otimes_k B^o)$).
Il n'est pas difficile de voir que pour un schéma $X$ de type fini et séparé sur $k$, 
le nc-schéma $\D(X)$ est propre si et seulement si $X$ est propre sur $k$, et de même 
$X$ est lisse sur $k$ si et seulement si le nc-schéma $\D(X)$ est lisse. Ces notions
de propreté et lissité sont donc compatibles avec les notions usuelles pour les schémas. \\

\textbf{Homologie de Hochschild des nc-schémas saturés.} 
Une des propriétés fondamentales des objets dualisables est qu'ils sont conservés par 
n'importe quel foncteur monoïdal symétrique. Par exemple, le foncteur d'homologie de Hochschild
$$C_* : k-\Sch^{nc} \longrightarrow L(k)^{op}$$
se promeut en un $\s$-foncteur monoïdal symétrique. En particulier, si $T$ 
est un nc-schéma saturé alors $C_*(T)$ est objet dualisable de $L(k)$, ou de manière
équivalente c'est un complexe parfait. En particulier, 
la classe $\sum_{i} (-1)^i[HH_i(T)]$ est bien définie dans $K_0(k)$. Il en est de même 
pour les versions à coefficients: si $T$ est saturé et $f : T \longrightarrow T$
est un endomorphisme de $T$, alors $C_*(T,f)$ est un complexe parfait. On dispose ainsi 
d'une classe de K-théorie $\sum_{i} (-1)^i[HH_i(T,f)]$. En particulier, lorsque $S=Spec\, k$ 
est connexe, on dispose d'une fonction rang $rang : K_0(k) \longrightarrow \mathbb{Z}$, et 
l'on peut alors définir la caractéristique d'Euler
$$\chi(HH_*(T,f)):=rang(\sum_{i} (-1)^i[HH_i(T,f)]) \in\mathbb{Z}.$$
Il est intéressant de noter que le complexe $C_*(T,f)$ peut aussi s'interpréter
comme la \emph{trace de f} dans l'$\s$-catégorie monoïdale symétrique $k-\Sch^{nc}$. 
De manière générale, si $x$ est un objet dualisable dans une catégorie monoïdale 
symétrique $(C,\otimes)$, et si $f$ est un endomorphisme de $x$, on peut définir
une trace de $f$ comme un élément du monoïde $End(\mathbf{1})$ des endomorphismes
de l'unité. Par dualité, l'endomorphisme $f$ définit un morphisme $\Gamma(f) : \mathbf{1} 
\longrightarrow x\otimes x^\vee$, appelé le graphe de $f$, que l'on compose
avec le morphisme évaluation $ev : x \otimes x^\vee \longrightarrow 
\mathbf{1}$ pour obtenir le trace de $f$ 
$$Tr(f):=ev \circ \Gamma(f) : \mathbf{1} \longrightarrow \mathbf{1}.$$
Dans le cas où $C=k-\Sch^{nc}$, l'unité est la dg-catégorie $k$ (un seul objet et $k$ 
comme endomorphisme), et l'on peut voir que $End(k)$ s'identifie à l'espace
des complexes parfaits de $k$-modules. Ainsi, la trace d'un endomorphisme $f$ 
d'un nc-schéma saturé $T$ est un complexe de $k$-module, et l'on vérifie que l'on a
$C_*(T,f) = Tr(f)$ en tant qu'endomorphisme de l'unité $k$. 

Signalons enfin ce qu'il se passe dans le cas particulier d'une variété propre et lisse
$X$ sur un corps $k$ et d'un endomorphisme $f$ de $X$. Tout d'abord, on peut voir que
$\D(X)$ est autodual, par l'équivalence $\D(X) \simeq \D(X)^o$ qui envoie un 
complexe parfait $E$ sur $X$ sur le complexe dual $\mathbb{R}\underline{Hom}(E,\OO_X)$. 
Avec cette identification, le morphisme d'évaluation $\D(X) \otimes_k \D(X) \longrightarrow 
k$, correspond au dg-foncteur $k \longrightarrow \D(X) \otimes_k \D(X) \simeq \D(X\times_k X)$
qui pointe le faisceau structural de la diagonale $\Delta_X$. Le morphisme de coévaluation
quand à lui correspond au dg-foncteur $\D(X\times_k X) \longrightarrow \D(k)$ qui envoie
un complexe $E \in \D(X\times X)$ sur $\mathbb{H}(X\times_k X, \Delta_X \otimes^{\mathbb{L}} E)
$. Le fait que ces morphismes fasse de $\D(X)$ un objet dual de lui-même, implique 
qu'il induit un accouplement non-dégénéré en homologie de Hochschild
$$<-,-> : HH_*(\D(X)) \otimes_k HH_*(\D(X)) \longrightarrow k.$$
A l'aide de l'identification HKR $HH_*(\D(X)) \simeq H^*(X,\Omega^*_X)$, l'accouplement ci-
dessus n'est autre que le produit d'intersection en cohomologie de Hodge. Le fait que ce
produit d'intersection soit non-dégénéré est bien la dualité de Poincaré en cohomologie de 
Hodge. \\

\textbf{Cohomologie $\ell$-adique des nc-schémas saturés.}  On s'attend à ce que la cohomologie
$\ell$-adique des nc-schémas saturés se comportent de manière similaire à celle
des variétés propres et lisses. Cependant, ce fait se résume aujourd'hui en un ensemble
de questions ouvertes, que nous n'appellerons pas conjectures par prudence, et il manque
encore des énoncés généraux. Plutôt que d'en faire une liste exhaustive, citons ci-dessous
la propriété clé pour la formule des traces. 

Pour cela, rappelons qu'un un nc-schéma $T$ nous associons une théorie cohomologique
généralisée $K^T$. Les produits externes en K-théorie induisent des produits
$$K^T \wedge_{\BU_k}K^{T'} \longrightarrow K^{T\otimes_k T'}.$$
Après réalisation $\ell$-adique on obtient des morphismes de Kunneth
$$\mathbb{H}(T,\Ql) \otimes_{\Ql(\beta)} \mathbb{H}(T',\Ql) \longrightarrow 
\mathbb{H}(T\otimes_k T',\Ql).$$

\begin{df}\label{d5}
Soit $T$ un schéma non-commutatif sur $k$. Nous dirons que $T$ est \emph{$\otimes$-admissible}
si le morphisme du Kunneth 
$$\mathbb{H}(T,\Ql) \otimes_{\Ql(\beta)} \mathbb{H}(T^o,\Ql) \longrightarrow 
\mathbb{H}(T\otimes_k T^o,\Ql)$$
est un quasi-isomorphisme dans $\D(S,\Ql)$.
\end{df}

Nous nous attendons en réalité à ce que tout nc-schéma saturé $T$ soit $\otimes$-admissible
au sens ci-dessous. Nous pensons plus généralement qu'une large classe de nc-schémas
soient $\otimes$-admissibles, comme par exemple les nc-schémas de \emph{type fini} au sens de 
\cite{tv}. Cependant, de tels énoncés nous paraissent inaccessibles pour l'instant. Noter que
si $T$ est saturé et $\otimes$-admissible, alors il s'en suit formellement que 
$\mathbb{H}(T,\Ql)$ est un $\Ql(\beta)$-module localement constant et de rang fini. 
Lorsque $k$ est un corps algébriquement clos cela est équivalent au fait que 
chaque $\mathbb{H}^i(T,\Ql)$ est de dimension fini. Il faut donc comprendre la notion
de $\otimes$-admissibilité comme une notion de finitude plutôt qu'une notion formelle 
liée aux structures monoïdales, elle est essentiellement équivalent à la finitude
de la cohomologie $\ell$-adique. 

Pour terminer, sans avoir accès à des énoncés généraux assurant la $\otimes$-admissibilité, il
nous faudra vérifier au cas pas cas que les nc-schémas aux quels l'on souhaite 
appliquer la formule des traces, sont $\otimes$-admissibles. \\

\textbf{Une formule des traces non-commutative I.} Nous arrivons enfin à l'énoncé de la
formule des traces dans le cadre non-commutatif. Pour cela, soit $T$ un nc-schéma
saturé et que nous supposerons $\otimes$-admissible au sens de la définition \ref{d5}. 
Soit $f : T \longrightarrow 
T$ un endomorphisme de $T$.
Nous supposons aussi 
que $S=Spec\, k$ est connexe, de sorte à ce que l'entier $\chi(HH_*(T,f))$ soit bien 
défini. Par ailleurs, la $\otimes$-admissibilité de $T$ assure que 
$\mathbb{H}(T,\Ql)$ est un $\Ql(\beta)$-module dualisable. Ainsi, l'endomorphisme $f$ 
induit sur $\mathbb{H}(T,\Ql)$ possède une trace bien définie qui est un élément 
de $H^0_{et}(S,\Ql(\beta))$. Si l'on suppose de plus que $S$ est strictement hensélien, 
on dispose d'un isomorphisme canonique $H^0_{et}(S,\Ql(\beta))\simeq \Ql$. 
La trace de $f$ est ainsi un élément bien défini
$Tr(f_{|\mathbb{H}(T,\Ql)}) \in \Ql$.

\begin{thm}[\cite{trace}]\label{t3}
Avec les hypothèses et notation ci-dessus, on a 
$$\chi(HH_*(T,f)) = Tr(f_{|\mathbb{H}(T,\Ql)}).$$
\end{thm}

\'Etonnamment le théorème \ref{t3} n'est pas difficile à démontrer et se déduit formellement 
des aspects monoïdaux et de dualité. Le point clé est l'existence d'une 
version non-commutative du caractère
de Chen $Ch : K_0(T) \longrightarrow \mathbb{H}^{0}(T,\Ql)$, et de ses propriété 
multiplicative et fonctorielle. L'existence de ce caractère de Chern est quand à lui 
une conséquences des travaux de fondements de Tabuada et Robalo sur les motifs non-commutatifs
(voir \cite{ro,ta}). Bien entendu, lorsque $T=\D(X)$ pour $X$ une variété propre et lisse 
sur un corps $k$, et $f$ est induite par un endomorphisme de $X$, la formule du théorème \ref{t3}
retrouve la formule des traces de Grothendieck-Lefschetz-Verdier (\cite[Exp. III]{glv}).
Nous renvoyons à \cite{trace} pour les détails sur la preuve de cette formule des traces. \\

\textbf{Une formule des traces non-commutative II.} La formule des traces ci-dessus n'est 
malheureusement pas suffisante pour nous. Il arrive en effet parfois que certains
nc-schémas $T$ ne soient pas saturés en tant que nc-schéma au-dessus de $k$, mais
uniquement au-dessus de certaines bases convenables. C'est une situation similaire
à un k-schéma $X$ qui ne serait pas propre et 
lisse au-dessus de $S=Spec\, k$, mais propre et lisse au-dessus
de $Y$ un second $k$-schéma. Un exemple typique est celui d'une dg-catégorie
$2$-périodique $T$, dont nous allons voir un exemple d'origine géométrique important
à la section suivante. Une telle dg-catégorie, considérée comme une nc-schéma sur $k$ n'est 
jamais propre (à moins qu'elle soit nulle) 
car ses complexes de morphismes étant $2$-périodiques ne sont jamais 
parfaits sur $k$. 

Nous sommes ainsi amenés à considérer la situation plus générale suivante. 
On se fixe un nc-schéma $B$ \emph{en comonoïdes}, c'est à dire que l'on suppose que $B$ vient
avec un morphisme diagonal $\Delta : B \longrightarrow B\otimes_k B$ dans
$k-\Sch^{nc}$. En termes de dg-catégories, cela signifie que $B$ est munie d'une structure
monoïdale (non-nécessairement symétrique) $\otimes : B\otimes_k B \longrightarrow B$. 
On peut alors définir les \emph{nc-schémas au-dessus de $B$} comme étant des nc-schémas
sur $k$ munis d'une coaction de $B$. En termes de dg-catégories il s'agit de 
dg-catégories $T$ munie d'un produit tensoriel externe
$B\otimes_k T \longrightarrow T$, faisant de $T$ un objet en $B$-modules (on parle aussi
de dg-catégories $B$-linéaires). On dispose ainsi qu'une $\s$-catégorie $B-\Sch^{nc}$ des
nc-schémas au-dessus de $B$. Lorsque $B=k'$ est une $k$-algèbre commutative, on 
retrouve bien entendu les nc-schémas au-dessus de $k'$, mais l'intérêt de
cette notion est précisément de s'autoriser des bases $B$ plus exotiques. 
Un exemple fondamental est le cas où $B=k[u,u^{-1}]$, avec $u$ une variable libre en degré $2$.
Dans ce cas les nc-schémas au-dessus de $B$ ne sont autre que les nc-schémas $2$-périodiques.
Un autre exemple est donné par le cas où $B$ est la catégorie des représentations
$k$-linéaires d'un groupe fini $G$, auquel cas les nc-schémas au-dessus de $B$ peuvent être vus 
comme des nc-schémas $G$-équivariants. 

On dispose donc d'une $\s$-catégorie $B-\Sch^{nc}$ des nc-schémas au-dessus de $B$. Il nous 
faut maintenant introduire une notion de saturation pour de tels objets et une nouvelle
complication entre en jeu. En effet, à moins que $B$ soit un comonoïde cocommutatif (c'est à 
dire que la dg-catégorie monoïdale sous-jacente soit monoïdale symétrique), il n'existe pas
de structure monoïdale naturelle sur $B-\Sch^{nc}$ (cela est à rapprocher de la situation
des modules sur un anneau non-commutatifs). En revanche, pour un nc-schéma $T$ au-dessus
de $B$, on peut donner un sens à $T\otimes_B T^o$, comme nc-schéma au-dessus de $k$. Cela 
permet de parler de morphismes évaluation et coévaluation pour $T \in B-\Sch^{nc}$, et ainsi
de définir la notion de \emph{nc-schémas saturés au-dessus de $B$}. Nous renvoyons à 
\cite{trace} pour les détails formels. 

Soit donc $T$ un nc-schéma saturé sur $B$. Les morphismes évaluation et coévaluation
définissent des morphismes de nc-schémas sur $k$
$$\xymatrix{
B \ar[r]^-{coev} & T\otimes_B T^o \ar[r]^-{ev} & k,}$$
dont le composé correspond à un objet de $B$. Cet objet est, par définition, 
$C_*(T/B)$, le complexe d'homologie de Hochchild de $T$ relatif à $B$. On dispose ainsi 
d'une classe bien définie
$$[C_*(T/B)] \in K_0(B),$$
à la quelle il faut penser comme étant $\sum_{i}(-1)^i[HH_i(T/B)]$, bien que cette
somme formelle n'ait plus vraiment de sens ici (car $HH_i(T/B)$ ne sont plus des
objets de $B$ en général). De la même manière, si $T$ est un muni d'un endomorphisme $f$, 
on dispose d'un version à coefficients $C_*(T/B,f)$ et de sa classe 
$[C_*(T/B,f)] \in K_0(B)$. On peut alors considérer le caractère de Chern de cette classe
pour obtenir
$$Ch([C_*(T/B,f)]) \in \mathbb{H}^0(B,\Ql).$$
Par ailleurs, on dispose d'une version relative à $B$ de la $\otimes$-admissibilité
au sens de la définition \ref{d5}. Un nc-schéma $T$ sur $B$ est dit $\otimes$-admissible si 
le morphisme de Kunneth
$$\mathbb{H}(T,\Ql)\otimes_{\mathbb{H}(B,\Ql)}\mathbb{H}(T^o,\Ql)
\longrightarrow \mathbb{H}(T\otimes_B T^o,\Ql)$$
est un quasi-isomorphisme. Pour un tel nc-schéma, $\mathbb{H}(T,\Ql)$
est un $\mathbb{H}(B,\Ql)$-module dualisable. La trace de l'endomorphisme
induit par $f$ possède donc un sens (voir \cite{trace})
$$Tr(f_{|\mathbb{H}(T,\Ql)}) \in \mathbb{H}^0(B,\Ql).$$
La formule des traces au-dessus de $B$ se lit alors comme suit.

\begin{thm}[\cite{trace}]\label{t3'}
Avec les hypothèses et notation ci-dessus, on a 
$$Ch([C_*(T/B,f)]) = Tr(f_{|\mathbb{H}(T,\Ql)}).$$
\end{thm}

Pour terminer, il existe aussi une extension de la formule ci-dessus 
lorsque maintenant l'endomorphisme $f$ de $T$ n'est pas plus
$B$-linéaire mais recouvre un endomorphisme $g$ de $B$ qui préserve
la structure de comonoïde. Dans ce cas la formule garde
un sens mais est une égalité dans $\mathbb{H}^{0}(C_*(B,g),\Ql)$, où 
$C_*(B,g)$ est le complexe de Hochschild de $B$ à coefficients dans $g$ (qui est 
lui-même une dg-algèbre). Nous verrons que disposer d'une formule pour un cadre
à ce point général est utile. 

\section{Application à la formule du conducteur}

Nous arrivons enfin à la dernière partie de ce texte dans la quelle nous allons
appliquer la formule des traces du Th\'eor\`eme 
\ref{t3'} à la formule du conducteur de Bloch. Nous allons expliquer
comme cette formule des traces implique la Conjecture \ref{cb} sous l'hypothèse additionnelle 
que la monodromie est unipotente. Nous indiquerons aussi comment nous pensons que 
le cas général en découle. \\

Nous revenons au contexte de la Conjecture \ref{cb}. Soit $A$ un anneau 
de valuation hensélien, et $X \longrightarrow S=Spec\, A$ un morphisme propre de schémas. 
On suppose que $X$ est un schéma régulier et que la fibre générique géométrique  $X_t$ est 
lisse sur $K^{sp}$. On choisit $\pi$ une uniformisante de $A$ et 
on note $k=A/(\pi)$ le corps résiduel qui est supposé parfait. 
La Conjecture \ref{cb} se ramène aisément
au cas où $k$ est algébrique clos, et nous supposerons donc que $k=k^{sp}$. 

Nous allons introduire un nc-schéma $MF(X,\pi)$, qui représente d'une 
certaine façon la "différence" $X_o-X_t$, entre la fibre spéciale et la fibre
générique, au quel nous appliquerons la formule des traces \ref{t3} 
pour $f=id$. Pour commencer, 
pour un schéma
$Y$ et une fonction $f \in \Gamma(Y,\OO_Y)$ sur $Y$, on définit un nc-schéma
$MF(Y,f)$ de la manière suivante. On commence par supposer que $Y=Spec\, B$ est un schéma
affine, et ainsi la fonction $f$ peut être vue comme un élément de $B$. On définit une
dg-catégorie $MF(B,f)$ dont les objets sont des quadruplets $E:=(E_0,E_1,u,v)$, où
\begin{itemize}

\item $E_0$ et $E_0$ sont des $B$-modules projectifs de type fini

\item $u : E_0 \longrightarrow E_1$ et $v : E_1 \longrightarrow E_0$ sont des morphismes
de $B$-modules

\item on a $uv=vu=\times f$. 
\end{itemize}

Pour deux tels objets $E$ et $E'$, on dispose d'un complexe 
de morphismes naturel $Hom(E,E')$ qui définissent de manière évidente
une dg-catégorie $MF(B,f)$. Lorsque $Y$ n'est plus nécessairement affine
on définit $MF(Y,f)$ par recollement (voir \cite{mf}). 

Dans notre situation, d'une famille de variétés algébriques $X \longrightarrow S$, 
l'uniformisante $\pi$ définit une fonction sur $S$, et par précomposition une fonction 
sur $X$ que nous notons encore $\pi$. Cela permet de donner un sens à la dg-catégorie
$MF(X,\pi)$ et au nc-schéma correspondant. A priori le nc-schéma est au-dessus
de $A$. Cependant, il est par définition $2$-périodique et donc ne peut être 
propre sur $A$. Il est certainement propre sur $A[u,u^{-1}]$ mais en revanche n'est pas
lisse sur $A[u,u^{-1}]$ car il est concentré sur le point fermé $Spec\, k \hookrightarrow S$. 
Pour faire de $MF(X,\pi)$ un nc-schéma saturé il nous faut introduire
une base relativement exotique que nous noterons $B$, construite de la manière suivante.

Le point fermé $x:=Spec\, k \longrightarrow S$ permet de définir un groupoïde
$G:=x\times_S x$ au-dessus de $x$. Il faut comprendre ici le produit fibré $x\times_S x$
au sens de la géométrie dérivée de \cite{ems} 
pour obtenir un groupoïde $G$ non-trivial sur $x$. L'objet
$G$ est donc un groupoïde dans les schémas dérivés qui opère sur le schéma $x$. Cette opération
est non-triviale et d'une certaine façon le quotient de $x$ par l'action de $G$ s'identifie
au complété formel $\hat{S}_x$ de $S$ le long de $x$ (voir \cite[Thm. 4.1]{ems} 
pour des résultats
plus précis allant dans ce sens). Le schéma dérivé $G$ est facile à décrire, on a
$G\simeq Spec\, k\oplus k[1],$
où $k\oplus k[1]$ est la $k$-algèbre simpliciale commutative qui est l'extension de carré
nul triviale de $k$ par $k[1]$, la suspension de $k$. En revanche, la structure de 
groupoïde au-dessus de $x$ est subtile, particulièrement lorsque $A$ est d'inégale 
caractéristique. En égale caractéristique $G$ est un groupe abélien dérivé qui opère 
trivialement sur $x$. En revanche, lorsque la caractéristique est mixte $G$ n'est plus
équivalent à un groupe, et l'on remarque déjà que les morphismes source et but
$G \rightrightarrows x$ ne sont pas égaux. 

L'unité de $G$ produit un point fermé $x \hookrightarrow G$, et ainsi un faisceau
gratte-ciel $k(x)$ sur $G$. On pose 
$$B^+:=\mathbb{R}\underline{Hom}(k(x),k(x)),$$
la dg-algèbre des endomorphismes (dérivés) du faisceau cohérent $k(x)$. Il est facile
de voir que $B^+$ s'identifie à $k[u]$ avec $u$ en degré $2$. Cependant, $B^+$, en tant
que nc-schéma, est aussi muni d'une structure de comonoïde non-triviale. Cette structure
est en réalité induite par une structure $E_2$ sur $B^+$, elle-même produite par le produit 
de convolution des cohérents sur $G$ (voir \cite{mf} pour plus de détails). On pose alors
$$B:=B^+[u^{-1}].$$
Le nc-schéma $B$ est naturellement un objet en comonoïde, et l'on montre le fait suivant. 

\begin{prop}\label{p1}
Le nc-schéma $MF(X,\pi)$ vit naturellement au-dessus de $B$. De plus, 
il est saturé en tant qu'objet de $B-\Sch^{nc}$. 
\end{prop}

Nous sommes maintenant en mesure de démontrer le théorème suivant, qui est 
un cas particulier nouveau de la Conjecture \ref{cb}. 

\begin{thm}\label{t4}
Le Conjecture \ref{cb} est vrai si on suppose de plus que la monodromie opère de manière
unipotente sur $H^*(X_t,\Ql)$.
\end{thm}

Signalons les grandes étapes de la preuve du théorème \ref{t4}. On commence par 
observer que $K_0(B)\simeq \mathbb{Z}$. Ainsi, la classe $Ch([C_*(MF(X,\pi)/B)]) \in K_0(B)$
peut s'identifier à un entier que nous noterons $\chi(HH_*(MF(X,\pi)/B))$. 

\begin{enumerate}

\item Avec les notations précédentes, et sans supposer que la monodromie est unipotente, 
on a 
$$\chi(HH_*(MF(X,\pi)/B)) = deg([\Delta_X.\Delta_X]_o).$$

\item Si l'action de la monodromie sur $H^*(X_t,\Ql)$ est unipotente, alors
$MF(X,\pi)$ est $\otimes$-admissible au-dessus de $B$ au sens de la définition \ref{d5}. 

\item Si l'action de la monodromie sur $H^*(X_t,\Ql)$ est unipotente, alors
on a 
$$\chi(\mathbb{H}^{*}(MF(X,\pi),\Ql))=\chi(X_o) - \chi(X_t).$$

\end{enumerate}

La formule des traces, associée aux trois faits précédents, implique alors que l'on a 
$$deg([\Delta_X.\Delta_X]_o) = \chi(X_o) - \chi(X_t).$$
Cette dernière formule est bien la Conjecture \ref{cb} car l'hypothèse d'unipotence assure
que le conducteur de Swan est nul. 

Il faut signaler qu'aucun des trois faits précédents n'est évident et chacun demande
une preuve à part entière. Par exemple, le lecteur pourra déduire le point (3) 
du résultat principal de \cite{mf}. \\

\textbf{Vers le cas général.} Pour terminer, signalons que nous pensons que la preuve que
nous donnons du théorème
\ref{t3} peut en réalité se généraliser pour donner une preuve de la Conjecture \ref{cb}. 
En effet, dans le cas général où la monodromie n'est plus unipotente, on sait qu'il
existe un revêtement fini $S' \longrightarrow S$, ramifié au point fermé $x$, 
tel que la famille $X':=X\times_S S' \longrightarrow S'$ soit à monodromie unipotente. 
Notons $G$ le groupe de Galois du revêtement $S' \rightarrow S$. Le groupe $G$ opère sur 
$X'$ et nous pouvons pour tout $g\in G$ appliquer la formule des traces \ref{t3'} (ou plutôt 
sa version généralisée où l'endomorphisme opère aussi non-trivialement sur la base $S'$).
Il se trouve que $X'$ n'est en général plus régulier, et le nc-schéma $MF(X',\pi')$ au quel
on applique la formule des traces demande alors à être défini par sa version cohérente
(voir \cite{ep}). 
Dans ce cas, la formule des traces ne calcule pas exactement ce que l'on voudrait, à 
savoir le caractère de l'action de $G$ sur la représentation virtuelle $H^*(X_0,\Ql)-H^*(X_t,
\Ql)$, mais plutôt une version homologique de cette représentation. C'est précisément 
dans la différence entre ces deux représentations que la caractère de Swan apparait. 
Cette approche est cependant actuellement en cours d'investigation.

\end{document}